\RequirePackage[l2tabu, orthodox]{nag} 

\documentclass[review,onefignum,onetabnum]{siamart190516}

\usepackage{amsmath} 
\usepackage{amssymb} 
\usepackage{amsfonts} 
\usepackage{amsopn} 

\usepackage{mathrsfs} 

\newsiamremark{remark}{Remark}
\newsiamremark{hypothesis}{Hypothesis}
\crefname{hypothesis}{Hypothesis}{Hypotheses}
\newsiamthm{claim}{Claim}

\usepackage{algorithm} 
\usepackage{algpseudocode} 

\usepackage{graphicx} 
\usepackage{epstopdf}
\ifpdf
  \DeclareGraphicsExtensions{.eps,.pdf,.png,.jpg}
\else
  \DeclareGraphicsExtensions{.eps}
\fi

\usepackage{threeparttable} 
\usepackage{booktabs} 

\usepackage[strict=true]{csquotes} 
\usepackage[strings]{underscore} 

\headers{Newton retraction}{Ruda Zhang}

\title{Newton retraction as approximate geodesics on submanifolds\thanks{Submitted to the editors.
    \funding{This work was supported by the National Science Foundation grant DMS-1638521.}
  }}

\author{Ruda Zhang\thanks{The Statistical and Applied Mathematical Sciences Institute, Durham, NC
    (\email{rzhang@samsi.info}).
    Department of Mathematics, North Carolina State University, Raleigh, NC
    (\email{rzhang27@ncsu.edu}).}
}

\ifpdf
\hypersetup{
  pdftitle={Newton retraction as approximate geodesics on submanifolds},
  pdfauthor={Ruda Zhang}
}
\fi

\nolinenumbers 

\begin{document}

\maketitle

\begin{abstract}
  Efficient approximation of geodesics is crucial for practical algorithms on manifolds.
  Here we introduce a class of retractions on submanifolds,
  induced by a foliation of the ambient manifold.
  They match the projective retraction to the third order
  and thus match the exponential map to the second order.
  In particular, we show that Newton retraction (NR) is always stabler
  than the popular approach known as oblique projection or orthographic retraction:
  per Kantorovich-type convergence theorems,
  the superlinear convergence regions of NR include those of the latter.
  We also show that NR always has a lower computational cost.
  The preferable properties of NR are useful for
  optimization, sampling, and many other statistical problems on manifolds.
\end{abstract}

  

\section{Introduction}
\label{sec:intro}

Consider a function $F \in C^k(\mathbb{R}^n, \mathbb{R}^c)$, for which 0 is a regular value.
By the regular level set theorem (see e.g. \cite[Thm 3.2]{Hirsch1976}),
the zero set $F^{-1}(0)$ is a properly embedded $C^k$ submanifold of $\mathbb{R}^n$
with codimension $c$ and dimension $d = n - c$.
We call $F(x)$ the constraint function and $\mathcal{M} = F^{-1}(0)$ the constraint manifold.
Depending on the context, $\mathcal{M}$ may also be called an implicitly-defined manifold,
a solution manifold, or an equilibrium manifold.
A geodesic $\gamma_v(t)$ in the manifold is a curve with zero intrinsic acceleration,
and all the geodesics are encoded in the exponential map $\exp: \mathscr{E} \mapsto \mathcal{M}$
such that $\exp(x, v) = \gamma_v(1)$,
where $(x, v) \in \mathscr{E} \subset T \mathcal{M}$ are the initial location and velocity.

The exponential map is crucial for analysis and computation on manifolds.
Application to problems on manifolds include
optimization \cite{Adler2002, Absil2008, GeR2015, ZhangHY2016b, Boumal2018},
differential equations \cite{Hairer2006},
interpolation \cite{Sander2015},
sampling \cite{Brubaker2012, Byrne2013, LiuC2016, Leimkuhler2016,
  Zappa2018, Mangoubi2018, Lelievre2019, Goyal2019, ZhangW2020},
and many other problems in statistics \cite{ChenYC2020}.
If the exponential map is not known in an analytic form or is computationally intractable,
it can be approximated by numerically integrating the geodesic trajectory,
i.e. solving the ordinary differential equation
$D_t \gamma_v'(t) = 0$ with initial conditions $\gamma_v(0) = x$, $\gamma_v'(0) = v$.
For submanifolds, this can also be done by projecting $x + v$ to $\mathcal{M}$,
which requires solving a constrained minimization problem.
In general, an approximation to the exponential map is called a retraction \cite{Adler2002}.
As is widely acknowledged (see e.g. \cite{Absil2008, ZhangHY2016b, Boumal2018}), 
retractions are often far less difficult to compute than the exponential map,
which allows for practical and efficient algorithms.

In this article, we introduce a class of second-order retractions on submanifolds,
which move points along manifolds orthogonal to the constraint manifold.
Of particular interest among this class is Newton retraction (NR),
which we show to have better convergence property and lower computational cost
than the popular approach known as oblique projection or orthographic retraction.
\Cref{tab:comparison} provides an overview of methods for computing geodesics,
where we summarize some of our results.

\begin{table*}[t]
  \label{tab:comparison}
  \centering
  \caption{Computation of geodesic steps on submanifolds, a qualitative comparison.}
  \begin{threeparttable}
    \begin{tabular}{llllll}
      \toprule
      method         & manifold   &  approx. & cost   &  stepsize & e.g.    \\
      \midrule
      analytic geodesics  & simple   & exact  & -    & any    & \cite{Byrne2013} \\
      numerical geodesics & level set & variable & high & tiny   & \cite{Leimkuhler2016} \\
      projective retraction & level set & 2nd    & high & almost any  &                   \\
      orthographic retraction & level set & 2nd    & low  & small  & \cite{Zappa2018} \\
      Newton retraction\tnote{*} & level set & 2nd (\ref{thm:rf-exp-2nd-order}) & lower (\ref{prop:faster}) & large (\ref{thm:stabler}) & \cite{Adler2002} \\
      \bottomrule
    \end{tabular}
    \begin{tablenotes}
    \item[*] Cross-referenced are results in this paper.
    \end{tablenotes}
  \end{threeparttable}
\end{table*}

\subsection{Related Literature}


Retraction is widely useful for optimization on manifolds.
Retraction was first defined in \cite{Adler2002}
for Newton's method for root finding on submanifolds of Euclidean spaces,
which was applied to a constrained optimization problem
over a configuration space of the human spine.
In particular, they proposed a retraction for constraint manifolds using Newton's method,
which we study in this paper.
Noisy stochastic gradient method \cite{GeR2015} escapes saddle points efficiently,
which uses projective retraction for constrained problems.
For geodesically convex optimization, \cite{ZhangHY2016b} studied several first-order methods
using the exponential map and a projection oracle, while acknowledging the value of retractions.
The Riemannian trust region method
has a sharp global rate of convergence to an approximate second-order critical point,
where any second-order retraction may be used \cite{Boumal2018}.


Sampling on manifolds often involves simulating a diffusion process,
which is usually done by a Markov Chain Monte Carlo (MCMC) method.
\cite{Brubaker2012} generalized Hamiltonian Monte Carlo (HMC) methods
to distributions on constraint manifolds,
where the Hamiltonian dynamics is simulated using RATTLE \cite{Andersen1983},
in which orthographic retraction is used to maintain state and momentum constraints.
\cite{Byrne2013} proposed geodesic Monte Carlo (gMC),
an HMC method for submanifolds with known geodesics.
\cite{LiuC2016} proposed two stochastic gradient MCMC methods for manifolds with known geodesics,
including a variant of gMC.
To sample on configuration manifolds of molecular systems,
\cite{Leimkuhler2016} proposed a scheme for constrained Langevin dynamics.
Criticizing the stability and accuracy limitations in the SHAKE method
and its RATTLE variant---both of which use orthographic retraction---their scheme
approximated geodesics by numerical integration. 
\cite{Zappa2018} proposed reversible Metropolis random walks on constraint manifolds,
which use orthographic retraction.
\cite{Lelievre2019} generalized the previous work to constrained generalized HMC,
allowing for gradient forces in proposal; it uses RATTLE.
In this line of work, it is necessary to check that the proposal is actually reversible,
because large timesteps can lead to bias in the invariant measure.
The authors pointed out that numerical efficiency in these algorithms
requires a balance between stepsize and the proportion of reversible proposed moves.
More recently, \cite{ZhangW2020} proposed a family of ergodic diffusion processes
for sampling on constraint manifolds, which use retractions defined by differential equations.
To sample the uniform distribution on compact Riemannian manifolds,
\cite{Mangoubi2018} proposed geodesic walk, which uses the exponential map.
\cite{Goyal2019} used a similar approach to sample the uniform distribution
on compact, convex subsets of Riemannian manifolds,
which can be adapted to sample an arbitrary density using a Metropolis filter;
both use the exponential map.

Many other statistical problems on constraint manifolds are discussed in \cite{ChenYC2020},
where Newton retraction can be very useful.
In probabilistic learning on manifolds \cite{ZhangRD2020EnvEcon},
a retraction based on constrained gradient flow is used
for inference and data augmentation on density ridge \cite{ZhangRD2020nbb}.


\section{Newton Retraction}
\label{sec:nr}

\subsection{Preliminaries}
\label{sub:prelim}

Retractions \cite{Adler2002}
are mappings that preserve the correct initial location and velocity of the geodesics;
they approximate the exponential map to the first order.
Second-order retractions \cite[Prop 5.33]{Absil2008} are retractions with zero initial acceleration;
they approximate the exponential map to the second order.
In general, we define retraction of an arbitrary order as follows.

\begin{definition}
  \label{def:retraction}
  \textit{Retraction} $R(x, v)$ of order $i$ on a $C^k$ manifold, $1 \le i < k$,
  is a $C^{k-1}$ mapping to the manifold from an open subset of the tangent bundle
  containing all the zero tangent vectors, such that at every zero tangent vector
  it agrees with the exponential map in Riemannian distance to the $i$-th order:
  $R \in C^{k-1}(U, M)$, $\zeta(M) \subset U \subset T \mathcal{M}$,
  $\forall (x, v) \in T \mathcal{M}$, $t \in \mathbb{R}$, $R(x, t v) = \exp(x, t v) + o(t^i)$,
  in the sense that, $\lim_{t \to 0} d_g(R(x, t v), \exp(x, t v)) / t^i = 0$.
\end{definition}

\cite{Absil2012} defined a class of retractions on submanifolds of Euclidean spaces,
called projection-like retractions.
This includes projective and orthographic retractions, both of which are second-order.

\begin{definition}[\cite{Absil2012}, Def 14]
  \label{def:projection-like}
  \textit{Retractor} $V(x, v)$ of a $C^k$ submanifold of $\mathbb{R}^n$ with codimension $c$
  is a $C^{k-1}$ mapping from tangent vectors to linear $c$-subspaces of the ambient space,
  such that for every zero tangent vector,
  affine space of the form $A(x, v) = x + v + V(x, v)$ intersects the submanifold transversely:
  $V \in C^{k-1}(U, G_{c, n})$, $\forall (x, v) \in U$, $A(x, v) \cap \mathcal{M} \ne \emptyset$,
  $\forall x \in \mathcal{M}$, $A(x, 0) \pitchfork \mathcal{M}$.
  Here $G_{c, n}$ is the Grassmann manifold. 
  \textit{Projection-like retraction} $R_V(x, v)$ induced by a retractor
  is a correspondence that takes a tangent vector to the set of points
  closest to the origin of the affine space that intersects the submanifold:
  $R_V: U \rightrightarrows \mathcal{M}$,
  $R_V(x, v) = \arg\min \{\|y - (x + v)\| : y \in A(x, v) \cap \mathcal{M}\}$.
  In particular, it is a mapping if the tangent vector is small enough:
  $\forall x \in \mathcal{M}$, $\exists U' \subset U$, $(x, 0) \in U'$:
  $R_V|_{U'} \in C^{k-1}(U', \mathcal{M})$.
\end{definition}

\begin{definition}
  \label{def:projective-orthographic}
  \textit{Projective retraction} $R_P(x, v) = P_{\mathcal{M}}(x+v)$,
  where projection $P_{\mathcal{M}}(x) = \arg\min\{\|y - x\| : y \in \mathcal{M}\}$,
  and can be seen as the projection-like retraction induced by retractor
  $V(x,v) = N_p \mathcal{M}$, where $p = P_{\mathcal{M}}(x+v)$.
  \textit{Orthographic retraction} $R_O(x, v)$ is the projection-like retraction
  induced by retractor $V(x, v) = N_x \mathcal{M}$.
\end{definition}

\begin{lemma}[\cite{Absil2012}, Thm 15, 22]
  \label{lem:projection-like}
  Projection-like retraction is a (first-order) retraction.
  It is second-order if the retractor maps every zero tangent to the normal space:
  $\forall x \in \mathcal{M}$, $V(x, 0) = N_x \mathcal{M}$.
\end{lemma}

\begin{figure}[t]
  \centering
  \includegraphics[width=.9\linewidth]{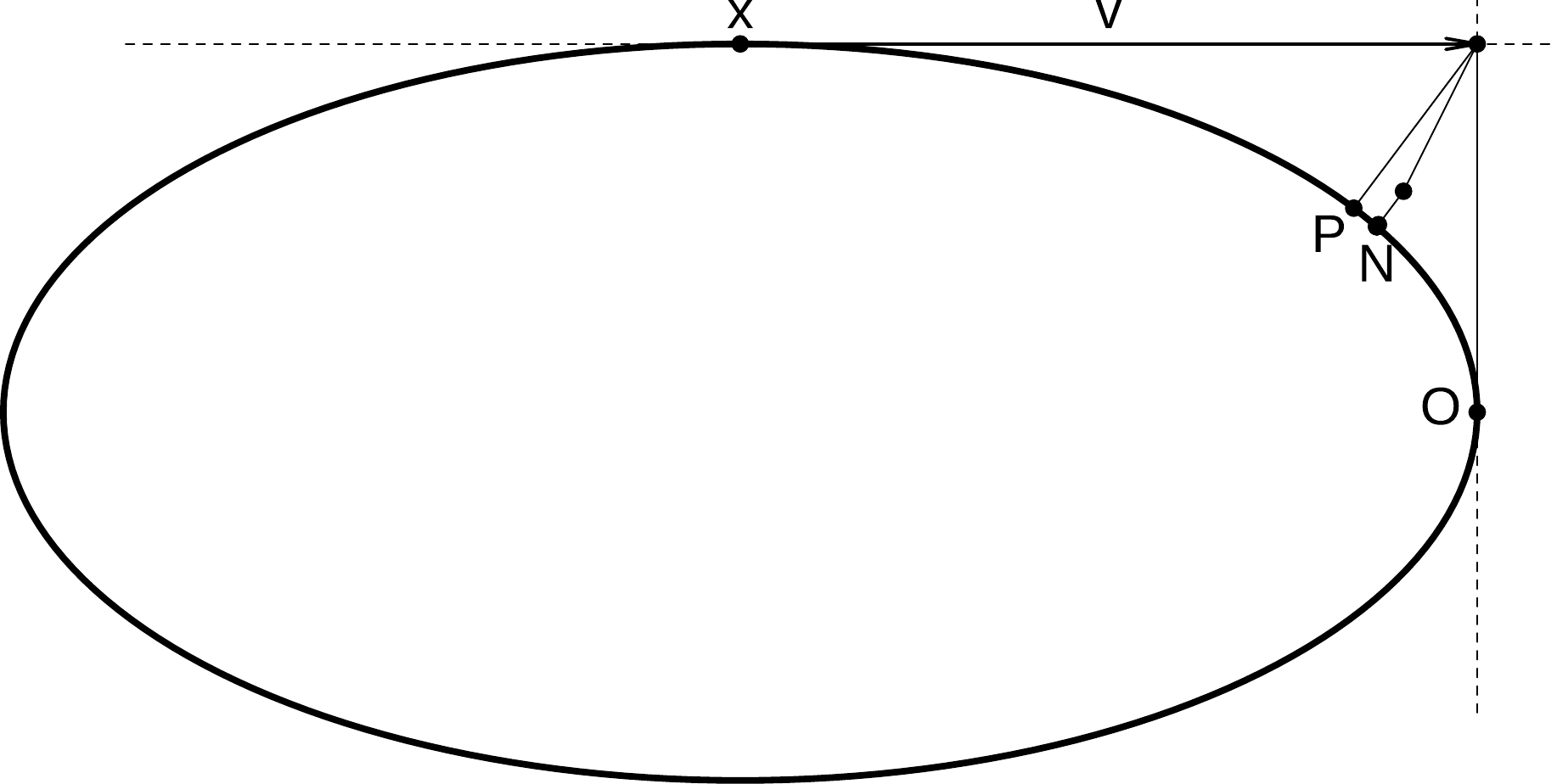}
  \caption{Illustration of retractions on an ellipse.
    With initial point $x$ and tangent vector $v$, retraction $R(x, v)$ returns:
    $O$, orthographic; $P$, projective; $N$, Newton (intermediate points are shown).}
  \label{fig:retractions}
\end{figure}

\subsection{Main Results}
\label{sub:main}

Here we define another class of second-order retractions on submanifolds,
based on foliations that intersect the submanifold orthogonally.
Such foliations may be generated by a dynamical system, 
where each leaf is a stable invariant manifold.
In particular, we are interested in Newton retraction,
generated by a discrete-time dynamical system with quadratic local convergence.
Such retractions can also be generated by flows:
\cite{ZhangW2020} used gradient flows of squared 2-norm of the constraint,
\cite{ZhangRD2020nbb} used constrained gradient flows of log density function;
both have linear local convergence and, with sufficiently small steps, global convergence.

\begin{definition}
  \label{def:retraction-foliation}
  \textit{Normal foliation} $\mathscr{F}$
  of a neighborhood of a codimension-$c$ submanifold of a Riemannian $n$-manifold
  is a partition of the neighborhood into connected $c$-submanifolds
  (called the \textit{leaves} of the foliation) which intersect the submanifold orthogonally:
  $\mathscr{F} = \sqcup_{p \in \mathcal{M}} \mathscr{F}_p$,
  $\cup_{p \in \mathcal{M}} \mathscr{F}_p = D$, $\mathcal{M} \subset D \subset \widetilde{\mathcal{M}}$,
  $\forall p \in \mathcal{M}$, $T_p \mathscr{F}_p = (T_p \mathcal{M})^\perp$.
  \textit{Retraction induced by a normal foliation}
  is the map $R_{\mathscr{F}} = \pi \circ \widetilde R$,
  where $\widetilde R$ is a retraction on the ambient manifold
  and $\pi: D \mapsto \mathcal{M}$ is the canonical projection:
  $\forall x \in D$, $x \in \mathscr{F}_{\pi(x)}$.
  If $\widetilde{\mathcal{M}} = \mathbb{R}^n$, let $\widetilde R$
  be the Euclidean exponential map $E(p, v) = p + v$,
  we have $R_{\mathscr{F}}(x, v) = \pi(x + v)$, $R_{\mathscr{F}}: E^{-1}(D) \mapsto \mathcal{M}$.
\end{definition}

Recall that for the under-determined system of nonlinear equations $F(x) = 0$,
Newton's minimum-norm step is $\delta(x) = - J^\dagger(x) F(x)$,
where Jacobian $J = \nabla F$ and $\dagger$ denotes the Moore-Penrose inverse.
If $J$ has full row rank, then $J^\dagger = J^T (J J^T)^{-1}$.
Newton map, or Newton's least-change update, is $N_F(x) = x + \delta(x) = x - J^\dagger(x) F(x)$.
Newton limit map is the mapping $N_F^\infty \in C^{k-1}(D, \mathcal{M})$ defined by
$N_F^\infty(x) = \lim_{n \to \infty} N_F^n(x)$, where $D$ is a neighborhood of $\mathcal{M}$.
\cite[Ex 4]{Adler2002} showed that given any retraction $\widetilde R$ on the ambient manifold,
the mapping $R = N_F^\infty \circ \widetilde R$ is a retraction on $\mathcal{M}$. 
We call this map Newton retraction. 

\begin{definition}
  \label{def:newton-retraction}
  \textit{Newton retraction}
  is the map $R_N \in C^{k-1}(E^{-1}(D), \mathcal{M})$ such that $R_N(x, v) = N_F^\infty(x + v)$,
  and can be seen as the retraction induced by the
  \textit{normal foliation determined by the Newton map}:
  $\mathscr{N} = \sqcup_{p \in \mathcal{M}} \mathscr{N}_p$,
  $\mathscr{N}_p = \{x \in D : N_F^\infty(x) = p\}$.
\end{definition}

To show that the retractions we defined are second-order, we give two lemmas.
First, projective retractions form a class of retractions of order two and not of any higher order.
Second, the retraction induced by a normal foliation
matches the projective retraction at zero tangent vectors to the third order.

\begin{lemma}
  \label{lem:rp-exp-2nd-order}
  For every $C^k$ submanifold $\mathcal{M}$, $k \ge 3$,
  the projective retraction $R_P$ satisfies:
  $\forall (x, v) \in T \mathcal{M}$, $t \in \mathbb{R}$, $R_P(x, t v) = \exp(x, t v) + o(t^2)$.
  There exists a $C^k$ submanifold, $k \ge 4$,
  such that the previous condition does not hold if $o(t^2)$ is replaced by $o(t^3)$.
\end{lemma}

\begin{lemma}
  \label{lem:rf-rp-3rd-order}
  Given a submanifold $\mathcal{M} \subset \mathbb{R}^n$ and
  a normal foliation $\mathscr{F}$ of a neighborhood of $\mathcal{M}$,
  $\forall (x, v) \in T \mathcal{M}$, $R_{\mathscr{F}}(x, v) = R_P(x, v) + o(\|v\|^3)$.
\end{lemma}

\begin{theorem}
  \label{thm:rf-exp-2nd-order}
  Retraction $R_{\mathscr{F}}$ induced by a normal foliation, which includes Newton retraction $R_N$,
  is a second-order retraction. This characterization of order is sharp.
\end{theorem}

Although $R_N$ and $R_O$ are both second-order retractions, they have different domain sizes.
Notice that the projection from a Euclidean space onto a compact subset
is uniquely defined except for a subset of measure zero,
so the projective retraction $R_P$ on a compact submanifold
is defined almost everywhere on the tangent bundle.
On the other hand, $R_O$ may be undefined for large tangent vectors
as the affine space $x + v + N_x \mathcal{M}$ fails to intersect the submanifold,
see \cref{fig:retractions}.
This precludes the domain of $R_O$ to a relatively small subset of the tangent bundle, 
regardless of implementation.
Since $R_N$ matches $R_P$ to the third order while $R_O$ and $R_P$ can differ on the third order,
it is easy to see that $R_N$ can have a larger domain than $R_O$.

\begin{algorithm}[t]
  \caption{Newton Retraction}
  \label{alg:nr}
  \begin{algorithmic}[1]
    \State Given: point and tangent vector $(x, v)$, convergence threshold $c_0$
    \State $x \gets x + v$
    \Repeat
    \State $J \gets J(x)$ \label{line:Jacobian}
    \State solve $(J J^T) y = F(x)$ \label{line:solve}
    \State $\delta \gets - J^T y$ \label{line:step}
    \State $x \gets x + \delta$
    \Until {$\|\delta\| < c_0$}
    \State \textbf{return} $x$
  \end{algorithmic}
\end{algorithm}

Now we characterize the domain of $R_N$ relative to that of $R_O$ in their usual implementation.
\Cref{alg:nr} gives an implementation of $R_N$,
which solves $F(\zeta) = 0$, $\zeta \in \mathbb{R}^n$ by Newton's method,
with initial value $x_0 = x + v$.
In comparison, $R_O$ is usually implemented by solving
$F(x_0 + J(x)^T y) = 0$, $y \in \mathbb{R}^c$ by Newton's method,
with initial value $y_0 = 0$.
This means replacing line~\ref{line:solve} with $(J J_{-1}^T) y = F(x)$
and line~\ref{line:step} with $\delta \gets - J_{-1}^T y$,
where $J_{-1} = J(x)$ is evaluated at the input $x$.
It can be seen as an augmented Jacobian algorithm
for solving under-determined systems of nonlinear equations:
denote the Stiefel manifold $V_{d, n} = \{X \in M_{n,d}(\mathbb{R}) : X^T X = I_d\}$,
given $V \in V_{d,n}$, step $\delta(x)$ is defined by $J(x) \delta(x) = - F(x)$, $V \delta(x) = 0$.
For $R_O$, the algorithm starts with $x_0 = x + v$ 
and $V$ satisfies $J(x) V = 0$. 
Kantorovich-type convergence theorems
for Newton's method and augmented Jacobian algorithms are given in \cite{Walker1990},
which provide sufficient conditions for immediately superlinear convergence.

\begin{definition}
  \label{def:hypotheses}
  Let $F \in C^1(\mathbb{R}^n, \mathbb{R}^m)$ and $J = \nabla F$.
  Let $C \subset \mathbb{R}^n$ be open and convex, $\alpha \in (0, 1]$, $K \ge 0$, and $B > 0$.
  We say function $F$ satisfies the normal flow hypothesis,
  $F \in \mathcal{H}_\text{NF}(C; \alpha, K, B)$,
  if $\forall x, y \in C$,
  (1) $\|J(x) - J(y)\| \le K \|x - y\|^\alpha$;
  (2) $\text{rank}(J(x)) = m$;
  (3) $\|J(x)^\dagger\| \le B$.
  Given $V \in V_{d,n}$, we say function $F$ satisfies the augmented Jacobian hypothesis,
  $F \in \mathcal{H}_\text{AJ}(V, C; \alpha, K, B)$ if $\forall x, y \in C$,
  (1) $\|J(x) - J(y)\| \le K \|x - y\|^\alpha$;
  (2') $\text{rank} \begin{bmatrix}J(x) \\ V\end{bmatrix} = n$;
  (3') $\left\| \begin{bmatrix}J(x) \\ V\end{bmatrix}^{-1} \right\| \le B$.
\end{definition}

\begin{theorem}[\cite{Walker1990}, Thm 2.1, 3.2]
  \label{thm:Walker1990}
  If $F \in \mathcal{H}_\text{NF}(C; \alpha, K, B)$
  then $\forall \eta > 0$, $\exists \epsilon > 0$, Newton's method satisfies:
  $\forall x_0 \in \{x \in C : B_\eta(x) \subset C, \|F(x)\| < \epsilon\}$,
  $\exists \zeta \in C \cap F^{-1}(0)$, $\lim_{k \to \infty} x_k = \zeta$;
  in particular, $\exists \beta > 0$, $\forall k \in \mathbb{N}$,
  $\|x_{k+1} - \zeta\| \le \beta \|x_k - \zeta\|^{1+\alpha}$.
  If $F \in \mathcal{H}_\text{AJ}(V, C; \alpha, K, B)$,
  then the previous statement holds for the augmented Jacobian algorithm.
\end{theorem}

Per the previous Kantorovich-type convergence theorem,
it follows immediately from the next lemma that
Newton retraction is always stabler than orthographic retraction.
Recall that $R_N$ has domain $E^{-1}(D)$, where $D$ is the domain of $N_F^\infty$,
i.e. the convergence region of Newton's method, and $R_O$ has domain $U$.
Let $U' \subset U$ be the convergence region of the usual implementation of $R_O$.

\begin{lemma}
  \label{lem:AJ-NF}
  For any $V \in V_{d,n}$, if $F \in \mathcal{H}_\text{AJ}(V, C; \alpha, K, B)$,
  then $F \in \mathcal{H}_\text{NF}(C; \alpha, K, B)$.
\end{lemma}

\begin{theorem}
  \label{thm:stabler}
  With the usual implementation of $R_N$ and $R_O$, for any $\alpha \in (0, 1]$,
  the order-$(1+\alpha)$ convergence region of $R_O$
  guaranteed by \cref{thm:Walker1990} is a subset of that of $R_N$:
  let $D_\alpha = \{x \in D : \exists (C, K, B, \eta, \epsilon),
  F \in \mathcal{H}_\text{NF}(C; \alpha, K, B), B_\eta(x) \subset C, \|F(x)\| < \epsilon\}$
  and $U'_\alpha = \{(x, v) \in U' : \exists (C, K, B, \eta, \epsilon),
  F \in \mathcal{H}_\text{AJ}(V_x, C; \alpha, K, B),
  B_\eta(x + v) \subset C, \|F(x + v)\| < \epsilon\}$, where $V_x \in V_{d,n}$, $J(x) V_x = 0$,
  then $\forall \alpha \in (0, 1]$, $U'_\alpha \subset E^{-1}(D_\alpha)$.
\end{theorem}

Our next result shows that Newton retraction is always faster than orthographic retraction.

\begin{proposition}
  \label{prop:faster}
  With the usual implementation of $R_N$ and $R_O$,
  for any $\alpha \in (0, 1]$, for any $(x, v) \in U'_\alpha$,
  the number of operations required for $R_N$ to converge is no greater than that for $R_O$.
\end{proposition}

\section{Discussion}
\label{sec:discussion}

The computational complexity per iteration in \cref{alg:nr} is dominated by
the evaluation of the Jacobian which are $c \times n$ real-valued functions,
and the linear solver in use, which invokes $O(c^a)$ algebraic operations, $a \approx 2.8$.
Since convergence is immediately superlinear (typically quadratic),
it is reasonable to assume that the algorithm stops after a fixed number of iterations.
Therefore the overall cost should suffice for most problems.

In case Jacobian evaluation is expensive and high numerical accuracy is unnecessary,
one may consider a modified Newton retraction (mNR):
run line~\ref{line:Jacobian} only for the first iteration,
denote the outcome as $J_0$, and replace line~\ref{line:solve} with $(J_0 J_0^T) y = F(x)$.
As a corollary of \cref{lem:projection-like}, mNR is a second-order retraction.
In this context, a natural implementation of $R_O$ is to use a chord method:
remove line~\ref{line:Jacobian}, and replace line~\ref{line:solve} with $(J_{-1} J_{-1}^T) y = F(x)$.
Both methods have linear convergence, but mNR has a faster rate:
let $\|x_{k+1} - x_\infty\| \le \mu \|x_k - x_\infty\|$,
then exists $\lambda \in (0, 1)$ such that $\mu = \lambda \|v\|^q$,
where $q = 1$ for $R_O$ and $q = 2$ for mNR, see e.g. \cite[Eq 6.2.15]{Allgower1990}.
By \cref{lem:AJ-NF}, mNR is no more stable than NR.


\section{Proofs}
\label{sec:proofs}

\begin{proof}[Proof of \cref{lem:rp-exp-2nd-order}]
  \cite[Ex 18]{Absil2012} showed that projective retractions are second-order retractions,
  so we only need to show that the projective retraction of some manifold is exactly second-order.
  Consider the circle $\mathbb{S}^1$ as a submanifold of the Euclidean plane,
  identified with the set $\{(\cos \theta, \sin \theta) : \theta \in [0, 2 \pi)\}$.
  Its exponential map is $\exp(x, v) = x e^{i \|v\|}$,
  and its projective retraction is $R_P(x, v) = (x + v) / \|x + v\|$.
  Without loss of generality, consider the point $x = (1, 0)$ and tangent vectors $v = (0, \theta)$.
  Now we can write the exponential map as $\exp(x, v) = e^{i \theta}$
  and the projective retraction as $R_P(x, v) = e^{i \arctan \theta}$.
  So the distance between them is $d(\exp, R_P) = 2 \sin ((\theta - \arctan \theta) / 2)$,
  which has a Taylor expansion at zero as $d(\exp, R_P) = \theta^3 / 3 + O(\theta^5)$.
  We can see that, for the unit circle, the projective retraction matches the exponential map
  up to the second order at zero tangent vectors.
\end{proof}

\begin{figure}[t]
  \centering
  \includegraphics[width=\linewidth]{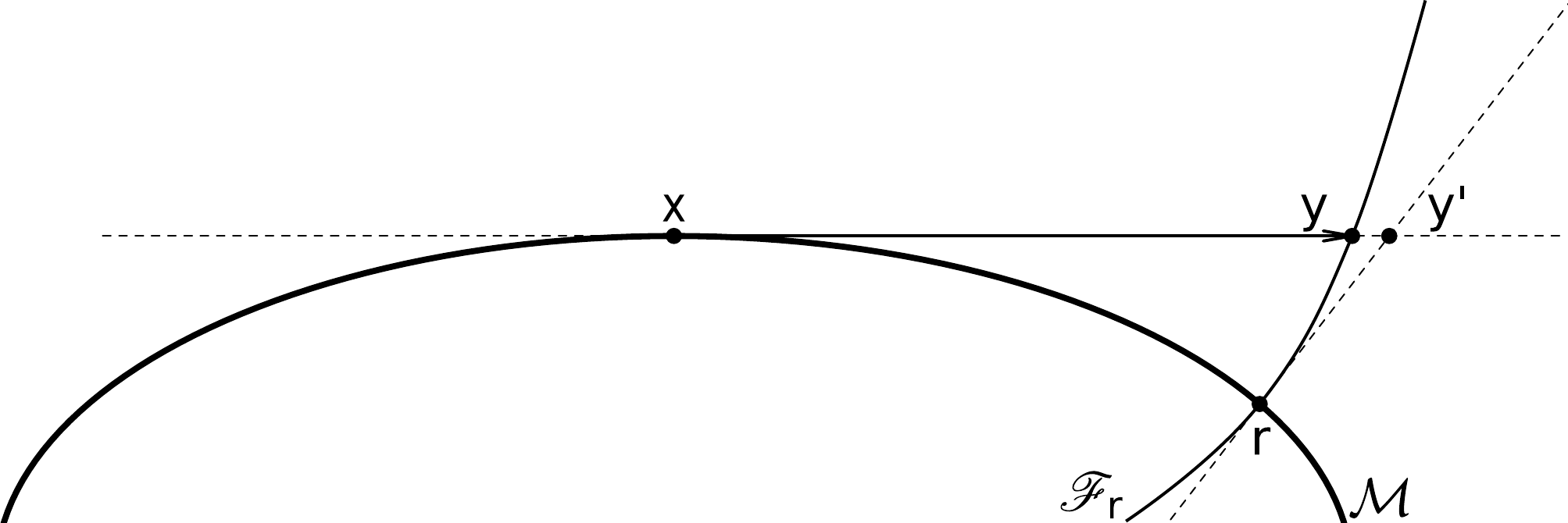}
  \caption{Approximate projection by a normal foliation.
    Leaf $\mathscr{F}_r$ intersects the submanifold $\mathcal{M}$ orthogonally at $N$
    and intersects tangent space $T_x \mathcal{M}$ at $y$.
    Tangent spaces $T_r \mathscr{F}_r$ and $T_x \mathcal{M}$ intersect at $y'$.}
  \label{fig:approximate-projection}
\end{figure}

\begin{proof}[Proof of \cref{lem:rf-rp-3rd-order}]
  For all $(x, v) \in T \mathcal{M}$ such that $y = x + v$
  is in the neighborhood of $\mathcal{M}$ that is partitioned by $\mathscr{F}$,
  define $r \in \mathcal{M}$ to be the unique point such that $y \in \mathscr{F}_r$.
  Note that $T_r \mathscr{F}_r = N_r \mathcal{M}$,
  so if $v = 0$ then $T_r \mathscr{F}_r$ and $T_x \mathcal{M}$ are orthogonal complements.
  Assume $v$ is small enough such that
  affine spaces $r + T_r \mathscr{F}_r$ and $x + T_x \mathcal{M}$ intersect transversely,
  define the unique point $y' = (r + T_r \mathscr{F}_r) \cap (x + T_x \mathcal{M})$.
  These constructs are illustrated in \cref{fig:approximate-projection}.
  Furthermore, define $v' = y' - x \in T_x \mathcal{M}$ and $u' = y' - r \in N_r \mathcal{M}$.
  Because $R_{\mathscr{F}}(x, v) = R_P(x, v')$, then $\exists w \in T_x \mathcal{M}$,
  $R_{\mathscr{F}}(x, v) = R_P(x, v) + \left(\frac{\partial R_P}{\partial v}(x, w)\right)(v' - v)$,
  that is, $R_{\mathscr{F}}(x, v) = R_P(x, v) + O(\|v' - v\|)$.
  To prove the theorem, we only need to show $\|v' - v\| = O(\|v\|^4)$.

  First we show that $\|u'\| = O(\|v'\|^2)$.
  Parameterize $\mathcal{M}$ at $x$ as the graph of a function $G: S_x \mapsto N_x \mathcal{M}$,
  $S_x \subset T_x \mathcal{M}$, such that $\forall v \in S_x$, $x + v + G(v) \in \mathcal{M}$.
  We see that $\|G(v')\| = O(\|v'\|^2)$.
  Let $p' = x + v' + G(v')$, because $d(y', \mathcal{M}) \le d(y', p')$, 
  we have $\|u'\| \le \|G(v')\|$.
  Thus, $\|u'\| = O(\|v'\|^2)$.

  Second, we show that $\|v' - v\| = O(\|u'\|^2)$.
  Parameterize $\mathscr{F}_r$ at $r$ as the graph of a function $L: S_r \mapsto N_r \mathscr{F}_r$,
  $S_r \subset T_r \mathscr{F}_r$, such that $\forall u \in S_r$, $r + u + L(u) \in \mathscr{F}_r$.
  We see that $\|L(u)\| = O(\|u\|^2)$.
  For all $v, w \in \mathbb{R}^n$, if $\|v\| \|w\| \ne 0$,
  define angle $\angle(v, w) = \arccos\left(\langle v, w \rangle / (\|v\| \|w\|)\right)$,
  otherwise define $\angle(v, w) = \pi / 2$.
  The first principal angle between linear subspaces $V, W \subset \mathbb{R}^n$ is defined as
  $\theta_1(V, W) = \min\{ \angle(v,w) : v \in V, w \in W \}$.
  Because $\theta_1(T_x \mathcal{M}, N_x \mathcal{M}) = \pi / 2$,
  we have $\theta_1(T_x \mathcal{M}, N_r \mathcal{M}) = \pi / 2 + O(\|v\|)$.
  Let $\beta = \angle(v - v', u - u')$,
  since $v - v' \in T_x \mathcal{M}$, $u - u' \in T_r \mathscr{F}_r = N_r \mathcal{M}$,
  we have $\beta \le \theta_1(T_x \mathcal{M}, N_r \mathcal{M}) = \pi / 2 + O(\|v\|)$.
  Thus, $\|v - v'\| = (\sin \beta)^{-1} \|L(u)\| = O(\|L(u)\|) = O(\|u\|^2)$
  and $\|u - u'\| = \|v - v'\| \cos \beta = \|v - v'\| O(\|v\|)$.
  Because $\|u\| \le \|u'\| + \|u' - u\|$,
  we have $\|u\| \le \|u'\| + o(\|v - v'\|) = \|u'\| + o(\|u\|^2)$,
  that is, $\|u\| = O(\|u'\|)$.
  We conclude that $\|v' - v\| = O(\|u'\|^2)$.

  Combining the previous two results gives $\|v' - v\| = O(\|v'\|^4)$.
  Because $\|v'\| \le \|v' - v\| + \|v\|$, we have $\|v'\| = O(\|v\|)$.
  This means $\|v' - v\| = O(\|v\|^4)$.
\end{proof}

\begin{proof}[Proof of \cref{thm:rf-exp-2nd-order}]
  Combining \cref{lem:rp-exp-2nd-order,lem:rf-rp-3rd-order},
  we have $R_{\mathscr{F}}(x, tv) = \exp(x, t v) + o(t^2) + o(t^3) = \exp(x, t v) + o(t^2)$,
  i.e. $R_{\mathscr{F}}$ is a second-order retraction.
  \cite[Thm 3.1]{Beyn1993} showed that $N_F^\infty$ induces a foliation of a neighborhood of
  $\mathcal{M}$ into $C^k$ $c$-submanifolds, which intersect $\mathcal{M}$ orthogonally.
  So $R_N$ fits \cref{def:retraction-foliation}
  as a retraction induced by a normal foliation, and thus it is second-order.
  Recall the circle example in the proof of \cref{lem:rp-exp-2nd-order},
  if $\mathbb{S}^1$ is defined as the zero set of $F(x) = \|x\|^2 - 1$, $x \in \mathbb{R}^2$,
  then $R_N = R_P$, which means in this case $R_N$ is only a second-order retraction.
\end{proof}

\begin{proof}[Proof of \cref{lem:AJ-NF}]
  By the definitions of the hypotheses,
  part (1) are identical,
  part (2) follows immediately from part (2'),
  so we only need to show part (3).
  For the rest of the proof, $x$ is an arbitrary point in $C$.
  To simplify notation, we will ignore explicit dependence on $x$,
  such that $J$ refers to $J(x)$, and so on.
  Since $J$ has full rank, we have QR decomposition $J^T = Q R$, where $Q \in V_{c,n}$ and
  $R \in U_+(c)$ is a upper triangular matrix of order $c$ with positive diagonal entries.
  Let $\bar{Q} = [Q, \tilde Q] \in O(n)$ be an orthogonal matrix of order $n$,
  whose first $c$ columns matches $Q$.
  Since $V \in V_{d,n}$, let $\bar{Q}_0 = [Q_0, \tilde Q_0] \in O(n)$, such that $V = \tilde{Q}_0^T$.

  First we show that $Q^T Q_0$ is non-singular and its spectral norm is no greater than 1.
  By (2'), $\begin{bmatrix}J \\ V\end{bmatrix} =
  \begin{bmatrix}R^T Q^T \\ \tilde Q_0^T \end{bmatrix}$ is non-singular.
  Since $R$ is invertible,
  this means $\begin{bmatrix} Q^T \\ \tilde Q_0^T \end{bmatrix}$ is non-singular,
  and thus $\begin{bmatrix} Q^T \\ \tilde Q_0^T \end{bmatrix} \bar{Q}_0 =
  \begin{bmatrix} Q^T \\ \tilde Q_0^T \end{bmatrix} [Q_0, \tilde Q_0]$ is non-singular.
  Note that $\tilde Q_0^T Q_0 = 0$, $\tilde Q_0^T \tilde Q_0 = I_d$,
  so $\begin{bmatrix} Q^T Q_0 & Q^T \tilde Q_0 \\ 0 & I_d \end{bmatrix}$ is non-singular,
  which means $Q^T Q_0$ is non-singular.
  Moreover, let $u \in \mathbb{S}^c$, then
  $\|Q^T Q_0 u\| \le \left\|\begin{bmatrix}Q^T \\ \tilde{Q}^T\end{bmatrix} Q_0 u\right\| =
  \|\bar{Q}^T Q_0 u \| = \| Q_0 u \| = 1$, which means $\rho(Q^T Q_0) \le 1$.

  Now we prove (3).
  By (3'), $\left\| \begin{bmatrix} J \\ V\end{bmatrix}^{-1} \right\| \le B$,
  that is, $\forall v \in \mathbb{R}^n$,
  $\left\| \begin{bmatrix}J \\ V\end{bmatrix}^{-1} v \right\| \le B \|v\|$.
  This means, $\forall w \in \mathbb{R}^n$,
  $\|w\| \le B \left\|\begin{bmatrix}J \\ V\end{bmatrix} w\right\|$.
  Equivalently, $\forall w \in \mathbb{S}^n$,
  $\left\|\begin{bmatrix}J \\ V\end{bmatrix} w\right\| \ge \frac{1}{B}$.
  Because $\forall u \in \mathbb{S}^c$, $Q_0 u \in \mathbb{S}^n$,
  so the previous inequality holds for $w = Q_0 u$.
  Note that $V Q_0 u = \tilde{Q}_0^T Q_0 u = 0$,
  the inequality becomes $\|J Q_0 u\| = \|R^T Q^T Q_0 u\| \ge \frac{1}{B}$.
  We have shown that $Q^T Q_0$ is non-singular and non-expansive,
  so $\forall z \in \mathbb{S}^c$, $\|R^T z\| \ge \frac{1}{B}$,
  or equivalently, $\forall u \in \mathbb{R}^c$, $\|R^T u\| \ge \frac{1}{B} \|u\|$.
  Define $\tilde u = R^{-1} u$, then $\forall \tilde u \in \mathbb{R}^c$,
  $\|R^T R \tilde{u}\| \ge \frac{1}{B} \|R \tilde{u}\|$.
  Define $\bar u = R^T R \tilde u$, then $\forall \bar u \in \mathbb{R}^c$,
  $\|R (R^T R)^{-1} \bar{u}\| \le B \|\bar{u}\|$.
  The left-hand side equals $\|Q R (R^T Q^T Q R)^{-1} \bar{u}\|$,
  that is $\|J^T(J J^T)^{-1} \bar{u}\|$ and in short $\|J^\dagger \bar{u}\|$.
  We conclude that $\|J^\dagger\| \le B$.
\end{proof}

\begin{proof}[Proof of \cref{thm:stabler}]
  $\forall (x, v) \in U'_\alpha$, $\exists (C, K, B, \eta, \epsilon)$:
  $F \in \mathcal{H}_\text{AJ}(V_x, C; \alpha, K, B)$.
  By \cref{lem:AJ-NF}, we have $F \in \mathcal{H}_\text{NF}(C; \alpha, K, B)$.
  Therefore, $x + v \in D_\alpha$ and $(x, v) \in E^{-1}(D_\alpha)$.
\end{proof}

\begin{proof}[Proof of \cref{prop:faster}]
  Since $U'_\alpha \subset E^{-1}(D_\alpha)$,
  the update sequences $\{x_k\}_{k \in \mathbb{N}}$ in $R_N$ and $R_O$ both satisfy
  $\|x_{k+1} - \zeta\| \le \beta \|x_k - \zeta\|^{1+\alpha}$.
  Because $d(x_0, \mathcal{M}) \le d(x_0, \mathcal{M} \cap (x_0 + N_x \mathcal{M}))$,
  the $x_k$ in $R_N$ will be remain closer to $\mathcal{M}$ than
  that in $R_O$ after the same number of iterations,
  and thus $R_N$ converges in no more iterations than $R_O$.
  Moreover, at each iteration, $R_N$ and $R_O$ both evaluate $F(x)$ and $J(x)$,
  and solve an order-$c$ system of linear equations, see line~\ref{line:solve}.
  But the coefficient matrix for $R_O$ is a generic matrix $J J_{-1}^T$,
  while that for $R_N$ is a symmetric matrix $J J^T$, which admits a faster linear solver.
  Therefore, the overall number of operations in $R_N$ is no greater than that in $R_O$.
\end{proof}



\bibliographystyle{siamplain}
\bibliography{nr.bib}

\end{document}